\documentclass[english,letterpaper,11pt,reqno]{amsart}

\usepackage[backend=bibtex,style=alphabetic,maxnames=6]{biblatex}
\usepackage{appendix}
\usepackage{amsbsy}
\usepackage{scalerel}
\usepackage{tcolorbox}
\usepackage{amsfonts}
\usepackage{amsmath}
\usepackage{amssymb}
\usepackage{amsthm}
\usepackage{graphicx}
\usepackage{ifthen}
\usepackage{textcomp}
\usepackage{enumitem, color, amssymb}
\usepackage{dsfont}
\usepackage{mathrsfs}
\usepackage{calrsfs}
\usepackage{accents}
\usepackage[bookmarksnumbered,colorlinks]{hyperref}
\hypersetup{colorlinks=true, linkcolor=blue, citecolor=blue}
\usepackage{cancel}
\usepackage{setspace}

\newcommand{\lra}{\longrightarrow}

\newcommand{\vep}{\varepsilon}

\makeatletter
\newcommand*{\defeq}{\mathrel{\rlap{%
                     \raisebox{0.25ex}{$\m@th\cdot$}}%
                     \raisebox{-0.25ex}{$\m@th\cdot$}}%
                     =}
\makeatother

\makeatletter
\newcommand*\owedge{\mathpalette\@owedge\relax}
\newcommand*\@owedge[1]{%
  \mathbin{%
    \ooalign{%
      $#1\m@th\bigcirc$\cr
      \hidewidth$#1\m@th\wedge$\hidewidth\cr
    }%
  }%
}
\makeatother

\newtheorem{thm}{Theorem}
\newtheorem{lemma}{Lemma}
\newtheorem{cor}{Corollary}

\newtheorem{defn}{Definition}
\newtheorem{prop}{Proposition}
\newtheorem*{definition-non}{Definition}
\newtheorem*{theorem-non}{Theorem}
\newtheorem*{proposition-non}{Proposition}
\newtheorem*{lemma-non}{Lemma}
\newtheorem*{corollary-non}{Corollary}

\newcommand{\beqa}{\begin{eqnarray}}
\newcommand{\beq}{\begin{equation}}
\newcommand{\eeqa}{\end{eqnarray}}
\newcommand{\eeq}{\end{equation}}

\newcommand\ipr[2]{\langle {#1},{#2}\rangle_{\!\scalebox{0.6}{$g$}}}

\newcommand\ww[2]{#1 \wedge #2}

\newcommand\imp{\hspace{.2in}\Rightarrow\hspace{.2in}}

\newcommand\comma{\hspace{.2in},\hspace{.2in}}
\newcommand\commas{\hspace{.1in},\hspace{.1in}}
\newcommand\commass{\hspace{.05in},\hspace{.05in}}

\newcommand\cp[1]{\hat{R}_{\scalebox{0.6}{$#1$}}}

\newtheorem{innercustomgeneric}{\customgenericname} %numbering theorems your own way
\providecommand{\customgenericname}{}
\newcommand{\newcustomtheorem}[2]{%
  \newenvironment{#1}[1]
  {%
   \renewcommand\customgenericname{#2}%
   \renewcommand\theinnercustomgeneric{##1}%
   \innercustomgeneric
  }
  {\endinnercustomgeneric}
}

\makeatletter
\def\Ddots{\mathinner{\mkern1mu\raise\p@
\vbox{\kern7\p@\hbox{.}}\mkern2mu
\raise4\p@\hbox{.}\mkern2mu\raise7\p@\hbox{.}\mkern1mu}}
\makeatother

\newcustomtheorem{customthm}{Theorem}
\newcustomtheorem{customlemma}{Lemma}

\newcommand{\CP}{\mathbb{CP}}

\newcommand{\dVg}{dV_{\scalebox{0.6}{$g$}}}

\newcommand{\R}{\mathbb R}
\newcommand{\G}{\mathsf G}
\newcommand{\ordsub}[1]{\scalebox{0.6}{$#1$}}
\newcommand{\Kt}[1]{K_{\ordsub{#1}}}
\newcommand{\Kht}[1]{\hat K_{\ordsub{#1}}}
\newcommand{\Rt}[1]{R_{\scalebox{0.6}{$#1$}}}
\newcommand{\Rht}[1]{\hat R_{\ordsub{#1}}}

\newcommand{\Kh}{\Kht{2n}}

\newcommand{\EK}{\mathcal E_{\scalebox{0.5}{$K$}}}
\newcommand{\KLK}{\mathcal{K}_{\scalebox{0.5}{$LK$}}}

\newcommand{\Bt}[1]{\mathcal B^{\scalebox{0.5}{$(#1)$}}}
\newcommand{\Bn}{\Bt{n}}
\newcommand{\dV}{\,dV_{\ordsub{g}}}
\newcommand{\tf}{\operatorname {tf}}

\newcommand{\tr}{\operatorname {tr}}
\newcommand{\Vol}{\operatorname {Vol}}
\newcommand{\Ric}{\operatorname {Ric}}
\newcommand{\scal}{\operatorname {scal}}
\newcommand{\id}{\operatorname {id}}

\begin{document}
\title{Higher Self-Duality}
\author{Amir Babak Aazami}
\address{Clark University\hfill\break\indent
Worcester, MA 01610}
\email{aaazami@clarku.edu}

\begin{abstract}
We use the totally trace-free component of Thorpe's middle curvature operator $\cp{2n}\colon\Lambda^{2n}\lra \Lambda^{2n}$ to define higher self-duality on oriented Riemannian $4n$-manifolds.  We show that higher self-dual metrics are absolute minimizers of a natural conformal energy, whose Bach tensor is symmetric, trace-free, and divergence-free, resolving a question posed by Kastor. We further show that the self- and anti-self-dual blocks give a signed-square formula for the top Pontryagin density. Compact self-dual examples include the Fubini--Study metric, though Einstein metrics need not be higher Bach-flat in general.
\end{abstract}

\maketitle

\section{Introduction}
\begingroup
\setlength{\parskip}{5pt}

On an oriented Riemannian 4-manifold $(M,g)$, the Hodge star $*$ is an involution on two-forms and therefore splits $\Lambda^2$ into a direct sum of its eigenspaces:
$
\Lambda^2=\Lambda^2_+\oplus\Lambda^2_-.
$
The Weyl curvature operator $\hat{W}\colon\Lambda^2\lra\Lambda^2$ preserves this splitting and thereby defines the restrictions $\hat W^{\pm}\colon \Lambda^2_\pm\lra \Lambda^2_\pm$.  A metric is called \emph{self-dual} when $\hat W^-=0$ and \emph{anti-self-dual} when $\hat W^+=0$; see, e.g., \cite{besse,lebrun2,taubes}.  But the blocks $\hat{W}^{\pm}$ are much more than just algebraic restrictions.  Indeed, the difference of their squared norms recovers the first Pontryagin number,
\beqa
\label{eqn:Hir}
p_1(M) = \frac{1}{4\pi^2}\int_M \left(|\hat{W}^+|^2-|\hat{W}^-|^2\right)\dVg,
\eeqa
while the vanishing of one block is exactly the integrability obstruction for the corresponding twistor space \cite{atiyahHitchinSinger}.  Thus $\hat{W}^{\pm}$ simultaneously connect curvature to both topology and complex geometry. There is also a variational story to $W$.  Namely, the \emph{Bach tensor} \cite{bach,besse} is the Euler--Lagrange tensor of the conformally invariant energy
$$
g\mapsto\int_M |W|^{2}\dV.
$$
It is symmetric, trace-free, divergence-free, and conformally covariant, and vanishes for both Einstein and self-dual/anti-self-dual metrics.  \mbox{It equals the} four-dimensional Fefferman--Graham obstruction tensor and is related to the variation of total $Q$-curvature \cite{feffermanGraham1985,feffermanGraham2012,grahamHirachi}. All in all, then, the geometry of the Weyl tensor in dimension four is particularly rich.

How much of this picture survives in dimensions $4n$?  To begin with, the middle-degree Hodge star $*$ still gives a splitting
$
\Lambda^{2n}=\Lambda^{2n}_+\oplus\Lambda^{2n}_-,
$
but $\hat{W}$ acts only on two-forms.  In this paper we use a construction in \cite{thorpe} to bridge this discrepancy.  Namely, for every even $p$, Thorpe alternated $p/2$ copies of the Riemann tensor to obtain a self-adjoint operator $\Rht{p}$ on $\Lambda^p$.  Kulkarni \cite{kulkarni} then isolated the totally trace-free part $\Kt{p}$ of $\Rt{p}$, which he called the $(p/2)^{\text{th}}$ conformal curvature tensor.  Here we use Thorpe's normalization and call $\Kt{2n}$ the \emph{Thorpe--Kulkarni tensor}. Its corresponding operator on $\Lambda^{2n}$ is denoted by $\Kh$, and its energy
\beqa
\label{eqn:EK}
\EK \defeq \int_M |\Kh|^{2}\dVg
\eeqa
is conformally invariant in all dimensions $4n$. In connection with Lovelock gravity, Kastor \cite{kastor} varied $\EK$, deriving the tensor
$$
\Bn_{ab}
=\Big((\Rt{2n-2})^{\mathbf i\mathbf j}
\nabla^\ell\nabla^k
-\frac{1}{2(2n-1)}
(\Rt{2n}^{[1]})^{\mathbf i k\mathbf j\ell}\Big)
(\Kt{2n})_{a\mathbf j\ell\,b\mathbf i k}.
$$
Here
$\mathbf i=(i_1,\ldots,i_{2n-2})$ and
$\mathbf j=(j_1,\ldots,j_{2n-2})$ are alternating multi-indices, $\Rt{0}=1$, and $\Rt{2n}^{[1]}$ denotes the first metric contraction of $\Rt{2n}$.  When $n=1$, this is the ordinary Bach tensor. However, for $n>1$ the question of whether $\Bn$ is symmetric, and hence is the full Euler--Lagrange tensor of \eqref{eqn:EK}, was posed in \cite{kastor}. We answer this question here in the affirmative.

\begin{theorem-non}
The tensor $\Bn$ is symmetric, and hence is the 
Euler--Lagrange tensor of
$
\EK(g)=\int_M|\Kh|^{2}\dV.
$
It is trace-free, divergence-free, and conformally covariant.
\end{theorem-non}
When $\Bn$ vanishes, we call the metric (\emph{higher}) \emph{Bach-flat}. We show that, in contrast to the four-dimensional case, Einstein metrics are not, in general, higher Bach-flat, but (higher) self-dual metrics are. Indeed, the total trace-freeness of $\Kh$ implies $[\Kh,*]=0$, and thus defines the blocks \mbox{$\Kh^\pm\colon\Lambda^{2n}\lra\Lambda^{2n}$;} this opens the door to a notion of higher self-duality.  We define a metric to be (\emph{higher}) \emph{self-dual} when $\Kh^-=0$ and \emph{anti-self-dual} when \mbox{$\Kh^+=0$.}  When $n=1$, $\Kt{2}=W$ and these definitions reduce to the classical ones. Not only are these metrics Bach-flat, but the blocks $\Kh^{\pm}$ always encode topological information: setting $
d_n\defeq\frac{(2n)!}{2^n n!}$ and $
\gamma_n\defeq(2\pi)^{2n}\left(\frac{2^n n!}{(2n)!}\right)^{\!2}
$, our main result is a direct generalization of \eqref{eqn:Hir}.

\begin{theorem-non}
On any closed oriented Riemannian $4n$-manifold $(M,g)$,
\beqa
\label{eqn:Hir2}
p_n(M)=\frac{d_n^2}{(2\pi)^{2n}}
\int_M\left(|\Kh^{+}|^{2}
-|\Kh^{-}|^{2}\right)\dVg.
\eeqa
If $\EK^{\pm}\defeq\int_M|\Kh^{\pm}|^{2}\,\dVg$, then
$$
\EK=\gamma_n|p_n(M)|+2\min\{\EK^{+},\EK^{-}\}.
$$
Hence $\EK\geq\gamma_n|p_n(M)|$, with equality if and only if the metric is higher self-dual or anti-self-dual.  Such metrics are absolutely minimizing and Bach-flat.
\end{theorem-non}

(The variation of $\int_M|\Rht{2n}|^2\dVg$ was recently considered in \cite{labbi25}.) For $n=1$, $p_1(M)=3\tau(M)$ and \eqref{eqn:Hir2} recovers Hirzebruch's four-dimensional signature formula \eqref{eqn:Hir}.  Note, however, that for $n>1$ the signature involves the full Hirzebruch $L$-polynomial (see \cite{besse}), whereas \eqref{eqn:Hir2} involves only the top Pontryagin number. Higher self-dual examples with $\Kt{2n}\neq 0$ exist, and include the Fubini--Study metric on $\mathbb{CP}^{2n}$.

The paper is organized as follows.  Section \ref{sec:operators} recalls Thorpe's definition of $\cp{2n}$ and Kulkarni's trace-free projection of it.  Section \ref{sec:bach} proves the symmetry and natural identities of the higher Bach tensor.  Section \ref{sec:self-duality} establishes automatic Hodge commutation and the Pontryagin identity \eqref{eqn:Hir2}.  Finally, Section \ref{sec:fubini-study} furnishes compact examples, while Section \ref{sec:products} exhibits a compact Einstein metric that is not higher Bach-flat, in contrast to the four-dimensional case.
\endgroup

\section{Higher curvature operators}
\label{sec:operators}

This section is a review of the material from \cite{thorpe,kulkarni,kastor,labbi2005} that will be needed to prove our main results; it is presented here for completeness and to fix notation and normalization constants. Let $(M,g)$ be a Riemannian $2m$-manifold, let $p\leq 2m$ be even, and let $R$ and $\hat{R}\colon \Lambda^2\lra \Lambda^2$ denote the Riemann curvature tensor and curvature operator, respectively. Thorpe's \emph{$p^{\text{th}}$ curvature tensor} is the covariant tensor defined by
\begin{align}
&\Rt{p}(u_1,\ldots,u_p,v_1,\ldots,v_p) \label{eq:Rp-def}\\
&\quad\defeq \frac{1}{2^{p/2}p!}
\sum_{\sigma,\tau\in S_p}\!\!\operatorname{sgn}(\sigma)
\operatorname{sgn}(\tau)
\prod_{r=1}^{p/2}
R(u_{\sigma(2r-1)},u_{\sigma(2r)},
  v_{\tau(2r-1)},v_{\tau(2r)}).\nonumber
\end{align}
The permutation $\sigma$ alternates only the $u$-variables, while $\tau$ alternates only the $v$-variables.
The factor $2^{p/2}$ compensates for
reversing the two entries in any curvature pair, and $p!$ compensates for
the repeated orderings of those pairs.  Thus \eqref{eq:Rp-def} is alternating in $u_1,\dots,u_p$, alternating in $v_1,\dots,v_p$, and symmetric under interchange of the
two blocks.  Its generalized algebraic and differential Bianchi
identities are, in indices,
\beqa
\label{eqn:gen_bian}
(\Rt{p})_{[i_1\ldots i_pj_1]j_2\ldots j_p}=0\comma
\nabla_{[i_0}(\Rt{p})_{i_1\ldots i_p]j_1\ldots j_p}=0.
\eeqa
It therefore defines a symmetric bilinear form on $\Lambda^p(TM)$, or equivalently a self-adjoint operator $\Rht{p}$, by
$$
\langle\Rht{p}(u_1\wedge\cdots\wedge u_p),
v_1\wedge\cdots\wedge v_p\rangle
\defeq\Rt{p}(u_1,\ldots,u_p,v_1,\ldots,v_p),
$$
where $\ipr{u_1\wedge \cdots \wedge u_p}{v_1\wedge \cdots \wedge v_p} \defeq \text{det}\big[g(u_i,v_j)\big]$. For $p=2$,
\begin{align*}
\Rt{2}(u_1,u_2,v_1,v_2)
={}&\frac14\bigl[
R(u_1,u_2,v_1,v_2)-R(u_2,u_1,v_1,v_2)\\
&\hspace{.55in}-R(u_1,u_2,v_2,v_1)+R(u_2,u_1,v_2,v_1)
\bigr]\\
={}&R(u_1,u_2,v_1,v_2),
\end{align*}
so that $\Rt{2}=R$ and $\Rht{2}=\hat R$.  The case $p=4$ is a homogeneous polynomial of degree $4/2=2$ in $R$; we record the full formula here because our eight-dimensional
example in Section \ref{sec:products} will use $\Rt{4}$ directly:
\begin{align}
&\Rt{4}(u_1,u_2,u_3,u_4,v_1,v_2,v_3,v_4) \label{eq:R4-explicit}\\
&=\frac13\big[
 R(u_1,u_2,v_1,v_2)R(u_3,u_4,v_3,v_4)
-R(u_1,u_2,v_1,v_3)R(u_3,u_4,v_2,v_4)\nonumber\\
&\quad+R(u_1,u_2,v_1,v_4)R(u_3,u_4,v_2,v_3)
+R(u_1,u_2,v_2,v_3)R(u_3,u_4,v_1,v_4)\nonumber\\
&\quad-R(u_1,u_2,v_2,v_4)R(u_3,u_4,v_1,v_3)
+R(u_1,u_2,v_3,v_4)R(u_3,u_4,v_1,v_2)\nonumber\\
&\quad-R(u_1,u_3,v_1,v_2)R(u_2,u_4,v_3,v_4)
+R(u_1,u_3,v_1,v_3)R(u_2,u_4,v_2,v_4)\nonumber\\
&\quad-R(u_1,u_3,v_1,v_4)R(u_2,u_4,v_2,v_3)
-R(u_1,u_3,v_2,v_3)R(u_2,u_4,v_1,v_4)\nonumber\\
&\quad+R(u_1,u_3,v_2,v_4)R(u_2,u_4,v_1,v_3)
-R(u_1,u_3,v_3,v_4)R(u_2,u_4,v_1,v_2)\nonumber\\
&\quad+R(u_1,u_4,v_1,v_2)R(u_2,u_3,v_3,v_4)
-R(u_1,u_4,v_1,v_3)R(u_2,u_3,v_2,v_4)\nonumber\\
&\quad+R(u_1,u_4,v_1,v_4)R(u_2,u_3,v_2,v_3)
+R(u_1,u_4,v_2,v_3)R(u_2,u_3,v_1,v_4)\nonumber\\
&\quad-R(u_1,u_4,v_2,v_4)R(u_2,u_3,v_1,v_3)
+R(u_1,u_4,v_3,v_4)R(u_2,u_3,v_1,v_2)\big].\nonumber
\end{align}

Aside from $\Rt{2}=R$, the top-degree case also explains Thorpe's normalization.  Indeed, if $p=2m$, then
$\Lambda^{2m}(T_xM)$ is one-dimensional, so $\Rht{2m}$ acts by a 
scalar.  This scalar is, up to sign, the Lipschitz--Killing curvature in
the Chern--Gauss--Bonnet theorem.  More precisely, let
$\{e_1,\ldots,e_{2m}\}$ be an oriented orthonormal basis, and let $\KLK$
denote the Lipschitz--Killing curvature, normalized to equal $1$ on
the unit round sphere.  Then
$$
\KLK
=(-1)^m\Rt{2m}(e_1,\ldots,e_{2m},e_1,\ldots,e_{2m}).
$$
The sign is present because $R(e_i,e_j,e_i,e_j)=-1$ on the unit round
sphere in this convention. Furthermore, the Pfaffian becomes
$$
\operatorname{Pf}(\Omega)
=\frac{2(2\pi)^m}{c_{2m}}
\Rt{2m}(e_1,\ldots,e_{2m},e_1,\ldots,e_{2m})\dV,
$$
where $c_{2m}$ is the volume of the $2m$-sphere. Thorpe's tensor still contains lower-order information obtained by
contracting one index from each alternating block.  Kulkarni 
removed all of this trace information, just as the Weyl tensor is obtained
by removing the Ricci and scalar-curvature parts of $R$.  Let us describe
this projection here.  Let
$\mathscr C_p(V)$ denote the space of tensors in
$\Lambda^pV^*\otimes\Lambda^pV^*$ which are symmetric under interchange
of the two factors and satisfy the generalized algebraic Bianchi identity
$$
T_{[i_1\ldots i_pj_1]j_2\ldots j_p}=0
$$
as in \eqref{eqn:gen_bian}. If
$T\in\mathscr C_p(V)$ and $\{e_1,\ldots,e_{2m}\}$ is orthonormal, define its
contraction $cT\in\mathscr C_{p-1}(V)$ by
$$
(cT)(u_1,\ldots,u_{p-1},v_1,\ldots,v_{p-1}) \defeq\sum_{i=1}^{2m}
T(e_i,u_1,\ldots,u_{p-1},e_i,v_1,\ldots,v_{p-1}),
$$
and write $T^{[r]}=c^rT$. This is independent of the chosen orthonormal basis; furthermore, because of the alternating properties of $p^{\text{th}}$ curvature tensors, $T$ is totally trace-free if and only if $cT=0$.  In the other direction, define metric insertion $\G T\in\mathscr C_{p+1}(V)$ by
\begin{align}
&(\G T)(u_1,\ldots,u_{p+1},v_1,\ldots,v_{p+1})\nonumber\\
&\qquad\defeq\sum_{i,j=1}^{p+1}(-1)^{i+j}g(u_i,v_j)
T(u_1,\ldots,\hat u_i,\ldots,u_{p+1},
  v_1,\ldots,\hat v_j,\ldots,v_{p+1}).\nonumber
\end{align}
For $p=1$, this is the usual
Kulkarni--Nomizu product:
$\G T = -g{\,\tiny \owedge\,} T$ (our sign convention for $\tiny \owedge$ is as in \cite{lee}).
More generally, Kulkarni showed that every
generalized curvature tensor is the orthogonal sum of a totally
trace-free part and pieces containing one or more explicit metric
insertions. Let us reproduce that result in our notation.

\begin{lemma}[\cite{kulkarni}]
\label{lem:projection}
Let $\dim V=2p$ and $T\in\mathscr C_p(V)$.  Its orthogonal projection onto the tensors annihilated by $c$ is
\begin{equation}
\label{eq:projection}
\tf(T)=\sum_{r=0}^{p}\frac{(-1)^r}{r!(r+1)!}\,\G^r c^rT.
\end{equation}
Moreover,
if $\mathscr K_q(V)\defeq
\ker\bigl(c\colon\mathscr C_q(V)\longrightarrow\mathscr C_{q-1}(V)\bigr)$ with $\mathscr K_0(V)=\R$,
then
\begin{equation}
\label{eq:orthogonal-decomp}
\mathscr C_p(V)=
\mathscr K_p(V) \oplus^\perp \G\mathscr K_{p-1}(V)
 \oplus^\perp\cdots\oplus^\perp \G^p\mathscr K_0(V).
\end{equation}
In particular, $\Kt{p}\defeq\tf(\Rt{p})$ is the top, totally trace-free summand of $\Rt{p}$.
\end{lemma}

\begin{proof}
For $A\in\mathscr C_q(V)$, consider $c(\G A)$; its terms separate into three types,
according to the locations of the contracted basis vector $e_i$.  When both 
copies of $e_i$ meet each other in $g$,
the sum over those terms gives $(\dim V)A$.  When one copy meets any of the
$q$ original indices $u_i$ in $g$, the sum over these terms gives $-qA$ (after using that \mbox{$\sum_ig(u_i,e_i)e_i=u_i$);} same for the $v_i$, for a total of $-2qA$.  In every remaining term the two $e_i$ enter $A$, and their sum is exactly $\G(cA)$.  Thus
$$
c(\G A)=\G(cA)+(\dim V-2q)A.
$$
Iterating this with $\dim V=2p$ and $T\in\mathscr C_p(V)$ gives
\begin{equation}
\label{eq:cGr}
c(\G^r c^rT)=\G^r c^{r+1}T+r(r+1)\G^{r-1}c^rT.
\end{equation}
The coefficient $r(r+1)$ can be seen by moving $c$
past the $r$ metric insertions $\G$ one at a time.  The first instance produces $(2p-2(p-1))\G^{r-1}c^rT$; the second produces an additional $(2p-2(p-2))\G^{r-1}c^rT$; the $j^{\text{th}}$ step produces an additional $(2p-2(p-j))\G^{r-1}c^rT$, and so on. Summing these gives
$$
\sum_{j=1}^r(2p-2(p-j))=2\sum_{j=1}^{r}j=r(r+1).
$$
Now put $a_r=(-1)^r/[r!(r+1)!]$.  Then
$a_r+(r+1)(r+2)a_{r+1}=0$, and an application of \eqref{eq:cGr} yields
\begin{align*}
c\left(\sum_{r=0}^{p}a_r\G^rc^rT\right)
={}&\sum_{r=0}^{p-1}a_r\G^rc^{r+1}T
+\sum_{r=1}^{p}a_rr(r+1)\G^{r-1}c^rT\\
={}&\sum_{r=0}^{p-1}
\bigl(a_r+(r+1)(r+2)a_{r+1}\bigr)
\G^rc^{r+1}T=0.
\end{align*}
This proves that the right-hand side of \eqref{eq:projection} lies in
$\ker c$; in particular, this is enough to ensure that $\Kt{p}\defeq\tf(\Rt{p})$ is totally trace-free. Now use that $\G$ and $c$ are $g$-adjoint,
$\langle\G A,B\rangle=\langle A,cB\rangle,$
hence $(\operatorname{im}\G)^\perp=\ker c$ (see \cite[Theorem~3.1]{labbi2005}).  In the sum
\eqref{eq:projection}, the $r=0$ term is $T$ itself and every term with
$r\geq1$ contains at least one factor of $\G$.  The difference between
$T$ and that sum therefore lies in $\operatorname{im}\G$.  Since the sum
itself lies in $\ker c=(\operatorname{im}\G)^\perp$, it is the
orthogonal projection.  Applying the same argument successively to the
lower-degree tensors gives \eqref{eq:orthogonal-decomp}.
\end{proof}

(The same contraction decomposition has been developed in the broader
calculus of curvature structures; see \cite{labbi2005}.)  Let us write out the cases $\Kt{2}$ and $\Kt{4}$. For $p=2$ and $\dim V=4$,
\beqa
\label{eqn:Ricci}
(cR)(Y,Z)=\sum_iR(e_i,Y,e_i,Z)=-\!\Ric(Y,Z)\comma c^2R=-\!\scal.
\eeqa
Then \eqref{eq:projection} yields
$$
\Kt{2}=R+\frac12\G(\Ric)-\frac1{12}\G^2(\scal)=R - \frac12 \Ric {\tiny \owedge}\,g + \frac{\scal}{12} g \,{\tiny \owedge}\, g = W,
$$
so that the first Kulkarni tensor is exactly the Weyl tensor. For $p=4$ and $\dim V=8$, however, \eqref{eq:projection} gives
$$
\Kt{4}=\Rt{4}-\frac12\G\Rt{4}^{[1]}
+\frac1{12}\G^2\Rt{4}^{[2]}
-\frac1{144}\G^3\Rt{4}^{[3]}
+\frac1{2880}\G^4\Rt{4}^{[4]}.
$$
Expanding the one-, two-, three-, and fourfold metric insertions with
weight-one antisymmetrization contributes the respective factors
$16,144,576,576$.  This changes the last four coefficients into
$-8,12,-4,1/5$, and gives
\begin{align}
(\Kt{4})_{abcd}{}^{efgh}
={}&(\Rt{4})_{abcd}{}^{efgh}
-8\delta_{[a}^{[e}(\Rt{4})_{bcd]i}{}^{fgh]i}+12\delta_{[a}^{[e}\delta_b^f(\Rt{4})_{cd]ij}{}^{gh]ij}\nonumber\\
&-4\delta_{[a}^{[e}\delta_b^f\delta_c^g(\Rt{4})_{d]ijk}{}^{h]ijk}+\frac15\delta_{[a}^{[e}\delta_b^f\delta_c^g\delta_{d]}^{h]}
(\Rt{4})_{ijkl}{}^{ijkl}.\label{eq:K4-index}
\end{align}
Here brackets have weight one: each full four-index antisymmetrization includes the factor $1/4!$. Examples aside, a key property is the following.
\begin{lemma}[\cite{kulkarni,kastor}]
For each even $p$ and \emph{$\text{dim}\,M=2p$},
$$
\int|\Kht{p}|^2\dV
$$
is conformally invariant.
\end{lemma}

\begin{proof}
Under $\tilde g=e^{2u}g$, we have $\tilde{\G}=e^{2u}\G$, $\tilde{c} = e^{-2u}c$, and the Riemann tensor has the form
$
\tilde R=e^{2u}(R+\G A_u),
$
where $A_u$ is a symmetric two-tensor formed from $\nabla^2u$,
$du\otimes du$, and $|du|^2g$.  Now form $\tilde R_{\ordsub{p}}$ by taking the $p/2$-fold product and alternating as in \eqref{eq:Rp-def}; this yields $\tilde R_{\ordsub{p}} = e^{pu}(R_{\ordsub{p}}+B_{\ordsub{p}})$, where $B_{\ordsub{p}}$ is a $p$-tensor (because both $\tilde R_{\ordsub{p}}$ and $R_{\ordsub{p}}$ are), comprised of $p/2$-fold products of the form $R \cdots R \cdot \G A_u\cdots \G A_u$ with at least one copy of $\G A_u$. Then
$$
\tilde K_{\ordsub{p}} = \tf(\tilde R_{\ordsub{p}}) = \tf\big(e^{pu}(R_{\ordsub{p}}+B_{\ordsub{p}})\big) = e^{pu}K_{\ordsub{p}} + e^{pu}\tf(B_{\ordsub{p}}).
$$
By construction, $\tf(B_{\ordsub{p}}) \in \ker c=(\operatorname{im}\G)^\perp$. However, every term in $\tf(B_{\ordsub{p}})$ is in $\operatorname{im}\G$, including $B_{\ordsub{p}}$ itself; indeed, $R^{p/2-1} \cdot \G A_u = \G(R^{p/2-1}\cdot A_u)$, and similarly for any $p/2$-fold product $R \cdots R \cdot \G A_u\cdots \G A_u$ with at least one $\G A_u$.
Thus $\tf(B_{\ordsub{p}})=0$ and $\tilde K_{\ordsub{p}}
=e^{pu}\Kt{p}$.  Raising the last $p$ indices uses $p$ copies of the
inverse metric, each of which contributes $e^{-2u}$.  Hence the associated
operator satisfies
$\Kht{p,\tilde g}=e^{-pu}\Kht{p,g}.$
In particular, $|\Kht{p}|^2$ has conformal weight
$e^{-2pu}$, whereas the volume form in dimension $2p$ has weight
$e^{2pu}$.
\end{proof}

\section{The higher Bach tensor}
\label{sec:bach}

Henceforth assume that $M$ is a closed oriented $4n$-manifold; set $p=2n$ and
\begin{equation}
\label{eq:energy}
\EK(g)\defeq\int_M|\Kh|^{2}\dV.
\end{equation}
We now examine the metric gradient $\Bn$ of this energy, defined for a smooth curve of metrics $g_t$ with $g_0=g$ and $\dot{g}_0=h$ by
$$
\left.\frac{d}{dt}\right|_{t=0}\EK(g_t)
=\beta_n\int_M\langle\Bn,h\rangle\dVg,
$$
where $\beta_n\neq 0$ is a universal constant. Its vanishing
will be the higher analogue of the Bach-flat equation. This tensor was derived in \cite{kastor} using a normalization in which his $n^{\text{th}}$ Riemann tensor is $2^n/(2n)!$ times Thorpe's $\Rt{2n}$.
We write it here using alternating multi-indices $\mathbf i=(i_1,\ldots,i_{2n-2})$ and
$\mathbf j=(j_1,\ldots,j_{2n-2})$; e.g., $a\mathbf j\ell$ means the ordered block
$a,j_1,\ldots,j_{2n-2},\ell$.  Repeated ordinary and multi-indices are
contracted using $g$.
Translating both curvature factors
and then discarding one harmless common constant gives
\begin{equation}
\label{eq:higher-bach}
\Bn_{ab}
=\Big((\Rt{2n-2})^{{\mathbf i}\mathbf j}
\nabla^\ell\nabla^k
-\frac{1}{2(2n-1)}
(\Rt{2n}^{[1]})^{{\mathbf i k}\mathbf j\ell}
\Big)(\Kt{2n})_{a\mathbf j\ell\,b\mathbf i k}.
\end{equation}
The two parts of this formula reflect the two parts of a curvature
variation.  Differentiating one of the $n$ curvature factors in
$\Rt{2n}$ leaves $n-1$ curvature factors, which combine to
$\Rt{2n-2}$; after two integrations by parts this produces the term with
two derivatives of $\Kt{2n}$.  Varying the metrics used to raise and
contract indices produces the algebraic term involving the first trace
$\Rt{2n}^{[1]}$.
For $n=1$, the multi-indices $\mathbf i,\mathbf j$ are empty,
$\Rt{0}=1$, $\Kt{2}=W$, and $\Rt{2}^{[1]}=cR=-\!\Ric$.  Thus
$$
\Bt{1}_{ab}
=\Big(\nabla^\ell\nabla^k
+\frac12\Ric^{k\ell}\Big)W_{a\ell bk},
$$
the ordinary Bach tensor in our curvature convention. Kastor showed that the symmetric part of \eqref{eq:higher-bach} is, up to
the universal normalization already discarded, the Euler--Lagrange tensor
of \eqref{eq:energy}.  The question of the symmetry of $\Bt{n}$, however, was left open for $n>1$
(see \cite[after Eq.~(36)]{kastor}).  We now settle this point.

\begin{thm}[Symmetry of the higher Bach tensor]
\label{thm:bach-symmetry}
The tensor $\Bn$ \mbox{in \eqref{eq:higher-bach}} is symmetric, and therefore is the Euler--Lagrange tensor of $\EK$. It is trace-free, divergence-free, and conformally covariant:
\begin{equation}
\label{eq:bach-properties}
g^{ab}\Bn_{ab}=0\comma
\nabla^a\Bn_{ab}=0\comma
\Bn_{\scalebox{0.6}{$e^{2u}g$}}
=e^{(2-4n)u}\Bn_{\scalebox{0.6}{$g$}}.
\end{equation}
\end{thm}

\begin{proof}
The second term of \eqref{eq:higher-bach} is algebraic, while the first
contains two covariant derivatives.  We prove their symmetry separately.
The only genuine difficulty is that interchanging the free indices in the
differential term reverses the order of those derivatives.  We shall show
that the curvature commutator measuring this reversal vanishes.

\emph{Step 1: the algebraic term.}
For the algebraic term,
the symmetry under interchange of the two alternating blocks in both
$\Rt{2n}^{[1]}$ and $\Kt{2n}$ gives
$$
\begin{aligned}
&(\Rt{2n}^{[1]})^{\mathbf i k \mathbf j\ell}
(\Kt{2n})_{b\mathbf j\ell\,a\mathbf i k}\\
&\qquad=(\Rt{2n}^{[1]})^{\mathbf i k \mathbf j\ell}
(\Kt{2n})_{a\mathbf i k\,b\mathbf j\ell}\\
&\qquad=(\Rt{2n}^{[1]})^{\mathbf j \ell \mathbf i k}
(\Kt{2n})_{a\mathbf j\ell\,b\mathbf i k}\\
&\qquad=(\Rt{2n}^{[1]})^{\mathbf i k\mathbf j\ell}
(\Kt{2n})_{a\mathbf j\ell\,b\mathbf i k},
\end{aligned}
$$
where in the second line we used the block symmetry of $\Kt{2n}$, in the third we relabeled $\mathbf i \leftrightarrow\mathbf j$ and $k \leftrightarrow\ell$, and in
the fourth we used the block symmetry of $\Rt{2n}^{[1]}$.  Thus the algebraic
 term is symmetric in $a,b$.

\emph{Step 2: rewrite the differential term as a double divergence.}
We would like to move both
derivatives in \eqref{eq:higher-bach} outside the contracted product. To that end, define
\begin{equation}
\label{eq:Pi}
\Pi_{a\ell bk}\defeq
(\Rt{2n-2})^{\mathbf i \mathbf j}
(\Kt{2n})_{a\mathbf j\ell\,b\mathbf i k}.
\end{equation}
Alternation and block interchange give
$$
\Pi_{a\ell bk}=-\Pi_{\ell abk}=-\Pi_{a\ell kb}\comma
\Pi_{a\ell bk}=\Pi_{bka\ell}.
$$
Thus $\Pi$ has the skew and pair-interchange symmetries of an ordinary
curvature tensor.  The differential Bianchi identity for $\Rt{2n-2}$ says
$$
\nabla_{[k}(\Rt{2n-2})_{i_1\ldots i_{2n-2}]}{}^{\mathbf j}=0\comma
\nabla^{[\ell}(\Rt{2n-2})_{\mathbf i}{}^{j_1\ldots j_{2n-2}]}=0.
$$
Because $\Kt{2n}$ is alternating in $(\mathbf i,k)$ and in
$(\mathbf j,\ell)$, contracting these identities with $\Kt{2n}$ and
$\nabla^k\Kt{2n}$ gives
$$
(\nabla^k(\Rt{2n-2})^{\mathbf i \mathbf j})
(\Kt{2n})_{a\mathbf j\ell\,b\mathbf i k}=0\hspace{.2in}\text{and}\hspace{.2in}(\nabla^\ell(\Rt{2n-2})^{\mathbf i \mathbf j})
\nabla^k(\Kt{2n})_{a\mathbf j\ell\,b\mathbf i k}=0.
$$
Indeed, because $\Kt{2n}$ alternates over
$i_1,\ldots,i_{2n-2},k$, contracting it against the first Bianchi sum
turns every summand of that identity into the same complete contraction,
with the sign already supplied by the alternation.  Since the Bianchi sum
is zero, that contraction is zero.  The second equality is the identical
argument in the other alternating block. Then
\begin{align}
\nabla^\ell\nabla^k\Pi_{a\ell bk}
={}&\nabla^\ell\bigl(
(\nabla^k\Rt{2n-2})^{\mathbf i \mathbf j}
(\Kt{2n})_{a\mathbf j\ell\,b\mathbf i k}\bigr)\nonumber\\
&\hspace{.1in}+\nabla^\ell\bigl(
(\Rt{2n-2})^{\mathbf i \mathbf j}
\nabla^k(\Kt{2n})_{a\mathbf j\ell\,b\mathbf i k}\bigr)\nonumber\\
={}&(\nabla^\ell(\Rt{2n-2})^{\mathbf i \mathbf j})
\nabla^k(\Kt{2n})_{a\mathbf j\ell\,b\mathbf i k}\nonumber\\
&\hspace{.1in}+(\Rt{2n-2})^{\mathbf i \mathbf j}
\nabla^\ell\nabla^k(\Kt{2n})_{a\mathbf j\ell\,b\mathbf i k}\nonumber\\
={}&(\Rt{2n-2})^{\mathbf i \mathbf j}
\nabla^\ell\nabla^k(\Kt{2n})_{a\mathbf j\ell\,b\mathbf i k}.\label{eq:double-div}
\end{align}
The first line vanishes before the outer derivative is taken, and the
first term in the second equation vanishes also.

\emph{Step 3: identify the tensor produced by commuting derivatives.}
We next prove the algebraic identity needed when the derivatives in
\eqref{eq:double-div} are commuted.  At a tangent space, consider the orthogonally invariant polynomial
$$
\Phi(R)\defeq|\tf(\Rt{2n})|^{2} = \ipr{\Kt{2n}}{\Kt{2n}}.
$$
If $S$ is another algebraic curvature tensor, let $(R+tS)_{\scalebox{0.6}{$2n$}}$ denote the $2n^{\text{th}}$ Thorpe constructed from $R+tS$ in place of $R$ in \eqref{eq:Rp-def}. Differentiating it yields
\beqa
\label{eq:last_step}
\left.\frac{d}{dt}\right|_{t=0}((R+tS)_{\scalebox{0.6}{$2n$}})_{a_1\ldots a_{2n}}{}^{b_1\ldots b_{2n}}
=\mu_n S_{[a_1a_2}{}^{[b_1b_2}
(\Rt{2n-2})_{a_3\ldots a_{2n}]}{}^{b_3\ldots b_{2n}]},
\eeqa
where $\mu_n\defeq n^2(2n-1)$ accounts for the normalization factors of $(R+tS)_{\scalebox{0.6}{$2n$}}$, $\Rt{2n-2}$ (recall \eqref{eq:Rp-def}), and those of both brackets; note that, after differentiating there will be $n$ terms, $SR\cdots R, RSR\cdots R,\dots,R\cdots RS$, but these are all equal after relabeling because both alternating blocks are fully antisymmetrized.  (The actual value of $\mu_n$ will not be needed below.) Next, using that $\tf((R+tS)_{\scalebox{0.6}{$2n$}})\big|_{t=0} = \tf(\Rt{2n}) = \Kt{2n}$, that $\frac{d}{dt}\big|_{t=0}\tf = \tf\frac{d}{dt}\big|_{t=0}$, that $\tf$ is self-adjoint, and finally that
$\tf(\Kt{2n})=\Kt{2n}$, we have
\begin{align}
(d\Phi)_{\scalebox{0.6}{$R$}}(S) = \left.\frac{d}{dt}\right|_{t=0}\Phi(R+tS)
&=2\left\langle
\left.\frac{d}{dt}\right|_{t=0}
\tf((R+tS)_{\scalebox{0.6}{$2n$}}),\Kt{2n}\right\rangle_{\!\!\!g}\nonumber\\
&=2\left\langle
\left.\tf\frac{d}{dt}\right|_{t=0}
(R+tS)_{\scalebox{0.6}{$2n$}},\Kt{2n}\right\rangle_{\!\!\!g}\nonumber\\
&=2\left\langle
\left.\frac{d}{dt}\right|_{t=0}
(R+tS)_{\scalebox{0.6}{$2n$}},\tf(\Kt{2n})\right\rangle_{\!\!\!g}\nonumber\\
&=2\left\langle
\left.\frac{d}{dt}\right|_{t=0}
(R+tS)_{\scalebox{0.6}{$2n$}},\Kt{2n}\right\rangle_{\!\!\!g}\nonumber\\
&\overset{\eqref{eq:last_step}}{=}2\mu_nS^{[a_1a_2[b_1b_2}(\Rt{2n-2})^{a_3\ldots a_{2n}]b_3\ldots b_{2n}]}\nonumber\\
&\hspace{.4in}\times (\Kt{2n})_{a_1\dots a_{2n}b_1\dots b_{2n}}\nonumber\\
&\overset{\eqref{eq:Pi}}{=}2\mu_nS^{abcd}\Pi_{abcd},\nonumber%\label{eq:Phi-derivative}
\end{align}
where in the last step, the anti-symmetrization brackets disappear when contracted against $\Kt{2n}$. Now define the two-tensor
$$
Q_{ab}\defeq\Pi_a{}^{cde}R_{bcde}.
$$
We show that $Q$ is symmetric.  The reason is orthogonal
invariance: an infinitesimal rotation cannot change $\Phi$.
Let $\alpha_{ab}$ be a skew-symmetric two-tensor, $\alpha_{ab}=-\alpha_{ba}$, and put
$A(t)=\exp(t\alpha)$. Note that $A(t) \in O(V)$, because $A(t)^TA(t) = \exp(t\alpha^*)\exp(t\alpha)=\exp(-t\alpha)\exp(t\alpha)=I$, where $\alpha = \alpha_a{}^b\colon V\lra V$ is skew-adjoint. Since $A(t)X = X+t\alpha X +O(t^2)$, differentiating the action
$(A(t)^* R)(X,Y,Z,W) = R(A(t)X,A(t)Y,A(t)Z,A(t)W)$ at $t=0$ yields
\beqa
\label{eqn:four}
\frac{d}{dt}\bigg|_{t=0}(A(t)^* R)=\alpha_a{}^eR_{ebcd}
+\alpha_b{}^eR_{aecd}
+\alpha_c{}^eR_{abed}
+\alpha_d{}^eR_{abce}.
\eeqa
Denote this by $(\alpha\cdot R)_{abcd}$. Then, because $\Phi$ is unchanged by orthonormal changes of frame,
$$
0=(d\Phi)_{\scalebox{0.6}{$R$}}(\alpha\cdot R) = 2\mu_n(\alpha\cdot R)^{abcd}\Pi_{abcd}.
$$
The four contractions arising on the right-hand side are equal; e.g.,
for the first two terms in \eqref{eqn:four}, skew-symmetry in $\Pi$ and $R$ gives
$$
\Pi^{abcd}\alpha_b{}^eR_{aecd}=-\Pi^{bacd}\alpha_b{}^eR_{aecd}=\Pi^{bacd}\alpha_b{}^eR_{eacd}
=\Pi^{abcd}\alpha_a{}^eR_{ebcd},
$$
where in the last step we interchanged the dummy indices $a,b$.  For the last two terms in \eqref{eqn:four}, pair-interchange symmetry in both $\Pi$ and $R$, together with another interchange of dummy indices, yields
\begin{align*}
\Pi^{abcd}\alpha_c{}^eR_{abed}
&=\Pi^{cdab}\alpha_c{}^eR_{edab}
=\Pi^{abcd}\alpha_a{}^eR_{ebcd},\\
\Pi^{abcd}\alpha_d{}^eR_{abce}
&=\Pi^{cdab}\alpha_d{}^eR_{ceab}
=\Pi^{dcab}\alpha_d{}^eR_{ecab}
=\Pi^{abcd}\alpha_a{}^eR_{ebcd}.
\end{align*}
Each of the four terms is therefore $\alpha^{ab}Q_{ab}$, and so
$
8\mu_n\alpha^{ab}Q_{ab}=0.
$
As $\alpha$ was arbitrary, the skew part of $Q$
must vanish: $Q_{ab}=Q_{ba}$.

\emph{Step 4: commute the derivatives.}
Pair-interchange symmetry of $\Pi$ and an interchange of dummy indices give
$$
\nabla^k\nabla^\ell\Pi_{a\ell bk}=\nabla^k\nabla^\ell\Pi_{bka\ell}
=\nabla^\ell\nabla^k\Pi_{b\ell ak}.
$$
It follows that
$$
[\nabla^\ell,\nabla^k]\Pi_{a\ell bk}=\nabla^\ell\nabla^k\Pi_{a\ell bk}
-\nabla^\ell\nabla^k\Pi_{b\ell ak}.
$$
With the convention
$([\nabla_{\!\partial_i},\nabla_{\!\partial_j}]\omega)_k=-R_{ijk}{}^m\omega_m$, the ordinary
Ricci commutation formula applied to all four covariant slots gives
\begin{align*}
[\nabla^\ell,\nabla^k]\Pi_{a\ell bk}
={}&-\!R^{\ell k}{}_{a}{}^m\Pi_{m\ell bk}-R^{\ell k}{}_{\ell}{}^m\Pi_{am bk}-R^{\ell k}{}_{b}{}^m\Pi_{a\ell mk}-R^{\ell k}{}_{k}{}^m\Pi_{a\ell bm}\\
={}&\frac12Q_{ba}+\Ric^{km}\Pi_{am bk}
-\frac12Q_{ab}-\Ric^{\ell m}\Pi_{a\ell bm}.
\end{align*}
Let us derive this carefully, working pointwise in an orthonormal frame, so that indices can be freely raised and lowered: $R^{\ell k}{}_{a}{}^m = R_{\ell k a m}$, etc.  For the first term, begin with the algebraic Bianchi identity with $a$ fixed, and then contract each term with $\Pi_{m\ell bk}$ (we suppress the summation notation in what follows):
$$
0=R_{\ell k a m}\Pi_{m\ell bk} + R_{k ma\ell}\Pi_{m\ell bk} + \!\!\!\!\underbrace{\,R_{m\ell a k}\Pi_{m\ell bk}\,}_{\text{$R_{a km\ell}\Pi_{bkm\ell} = Q_{ba}$}}.
$$
In the first term, interchange the dummy indices $m,\ell$, to obtain $R_{m k a \ell}\Pi_{\ell mbk}$; this equals the second term $R_{k ma\ell}\Pi_{m\ell bk} = -R_{mka\ell}\Pi_{m\ell bk} = R_{mka\ell}\Pi_{\ell m bk}$. Thus
$$
Q_{ba} = -2R_{\ell k a m}\Pi_{m\ell bk} = -2R^{\ell k}{}_{a}{}^m\Pi_{m\ell bk}.
$$
The third term follows in exactly the same way: $Q_{ab} = 2R^{\ell k}{}_{b}{}^m\Pi_{a\ell mk}$.
As $Q_{ab}=Q_{ba}$, these terms cancel. The two Ricci terms follow from $cR=-\!\Ric$. But a dummy-index interchange gives
$$
\Ric^{km}\Pi_{am bk}=\Ric^{km}\Pi_{ak bm}\comma
-\!\Ric^{\ell m}\Pi_{a\ell bm}=-\!\Ric^{km}\Pi_{ak bm},
$$
so they cancel also.  Thus $[\nabla^\ell,\nabla^k]\Pi_{a\ell bk}=0$. By \eqref{eq:double-div}, the differential term in \eqref{eq:higher-bach} is therefore symmetric. Thus $\Bn$ is symmetric.

\emph{Step 5: derive the three natural identities \eqref{eq:bach-properties}.} Since Kastor's variation
gives a nonzero universal multiple of $\Bn_{(ab)}$ and we have just
proved symmetry, $\Bn$ itself is a nonzero multiple of the full metric
gradient of $\EK$.  Equivalently, for one universal constant
$\beta_n\neq0$, our first-variation convention is
\begin{equation}
\label{eq:first-variation-convention}
\left.\frac{d}{dt}\right|_{t=0}\EK(g_t)
=\beta_n\int_M\langle\Bn,h\rangle\dV,
\end{equation}
where $g_t$ is a smooth curve of metrics satisfying $g_0=g$ and $\dot{g}_0=h$; $\beta_n$ is irrelevant and will cancel from every argument below.  For any
compactly supported function $\varphi$,
the conformal path $g_t=e^{2t\varphi}g = (1+t(2\varphi)+\mathcal{O}(t^2))g$ has $\dot g_0=h=2\varphi g$.  Conformal invariance of \eqref{eq:energy}
therefore gives
$$
0=2\beta_n\int_M\varphi\,g^{ab}\Bn_{ab}\dV,
$$
hence $g^{ab}\Bn_{ab}=0$ by the fundamental lemma of the calculus of variations.  Next, if $X$
is compactly supported and $g_t = \Phi_t^*g$ is the pullback of $g$ by its flow $\Phi_t$, then
$\dot g_0=\frac{d}{dt}\big|_{t=0}\Phi_t^*g=\mathcal L_Xg$ and diffeomorphism invariance gives
$$
0=2\beta_n\int_M(\Bn)^{ab}\nabla_aX_b\dV
=-2\beta_n\int_M\bigl(\nabla_a(\Bn)^{ab}\bigr)X_b\dV
$$
(recall that $M$ is closed). Hence $\nabla^a\Bn_{ab}=0$.  Finally, put $\tilde g=e^{2u}g$ and let
$h$ be an arbitrary compactly supported symmetric two-tensor.  Conformal
invariance applied to the paths $g+th$ and
$e^{2u}(g+th)=\tilde g+te^{2u}h$ gives
$$
\int_M\langle\Bn_{\scalebox{0.6}{$\tilde g$}},e^{2u}h\rangle_{\scalebox{0.6}{$\tilde g$}}
\,dV_{\!\scalebox{0.6}{$\tilde g$}}
=\int_M\langle\Bn_{\scalebox{0.6}{$g$}},h\rangle_{\!\scalebox{0.6}{$g$}}\,\dVg.
$$
Since $dV_{\scalebox{0.6}{$\tilde g$}}=e^{4nu}\dVg$
and two inverse metrics enter the inner product of covariant two-tensors,
we have
$
\langle\Bn_{\scalebox{0.6}{$\tilde g$}},e^{2u}h\rangle_{\scalebox{0.6}{$\tilde g$}}
=e^{-2u}\langle\Bn_{\scalebox{0.6}{$\tilde g$}},h\rangle_{\!\scalebox{0.6}{$g$}}.
$
Thus the integrand on the left is
$e^{(4n-2)u}\langle\Bn_{\scalebox{0.6}{$\tilde g$}},h\rangle_{\!\scalebox{0.6}{$g$}}$.
That $h$ was arbitrary therefore gives
$\Bn_{\scalebox{0.6}{$\tilde g$}}
=e^{(2-4n)u}\Bn_{\scalebox{0.6}{$g$}}$, and completes the proof.
\end{proof}

\section{Higher self-duality and the top Pontryagin form}
\label{sec:self-duality}

In dimension four, self-duality is possible because $\hat{W}$ preserves the two eigenspaces of the Hodge star on two-forms.
We now show that the middle Thorpe--Kulkarni operator $\Kht{2n}$ does the same
in all dimensions $4n$.  To begin with, $*\colon\Lambda^{2n}\lra \Lambda^{2n}$ satisfies $*^2=1$ and thus induces the splitting
$$
\Lambda^{2n}=\Lambda^{2n}_{+}
\oplus\Lambda^{2n}_{-}\comma
\Lambda^{2n}_{\pm}
\defeq\{\alpha:*\alpha=\pm\alpha\}.
$$
An arbitrary self-adjoint operator on $\Lambda^{2n}$ need not preserve
this splitting.  But we now show that total trace-freeness ensures
that the Thorpe--Kulkarni operator will do so.

\begin{lemma}[Automatic Hodge commutation]
\label{lem:hodge-commutation}
Let $T\in\mathscr C_{2n}(V)$ be totally trace-free on an oriented Euclidean
space of dimension $4n$.  Then $[\hat{T},*]=0$. 
In particular, $[\Kht{2n},*]=0$.
\end{lemma}

\begin{proof}
We show that $*\hat T*=\hat T$. Choose an oriented orthonormal basis, and let
$\mathbf i,\mathbf j,\mathbf k,\boldsymbol\ell$ denote ordered
$2n$-tuples.  Repeated bold lowercase indices are summed over all such
ordered tuples. If $e_{\mathbf i} \defeq e_{i_1}\wedge \cdots \wedge e_{i_{2n}}$ with $i_1< \cdots < i_{2n}$, then
$$
*\,e_{\mathbf i} = \frac{1}{(2n)!}\vep_{\mathbf i}{}^{{\mathbf j}}e_{{\mathbf j}} \comma \hat{T} e_{\mathbf i} = \frac{1}{(2n)!}T_{\mathbf i}{}^{{\mathbf j}}e_{{\mathbf j}}\comma (*\hat{T}*)e_{\mathbf i} = \frac{1}{(2n)!}(*\hat{T}*)_{\mathbf i}{}^{{\mathbf j}}e_{{\mathbf j}},
$$
where $\vep_{\mathbf i \mathbf j}$ is the Levi--Civita $4n$-symbol. As a consequence,
$$
(*\hat T*)e_{\mathbf i} = \frac{1}{((2n)!)^3}\vep_{\mathbf i}{}^{{\mathbf j}}T_{\mathbf j}{}^{{\mathbf k}}\vep_{\mathbf k}{}^{{\boldsymbol\ell}}e_{{\boldsymbol \ell}} \imp
(*\hat T*)_{\mathbf i}{}^{\boldsymbol \ell}
=\frac{1}{((2n)!)^2}\vep_{\mathbf i}{}^{{\mathbf j}}T_{\mathbf j}{}^{{\mathbf k}}\vep_{\mathbf k}{}^{{\boldsymbol\ell}}.
$$
Expanding the sum on the right-hand side yields
$$
(*\hat{T}*)_{\mathbf i\boldsymbol\ell}
 =a\,T_{\mathbf i^c\boldsymbol\ell^c}\comma a\defeq\varepsilon_{\mathbf i\mathbf i^c}
   \varepsilon_{\boldsymbol\ell\boldsymbol\ell^c},
$$
so it remains to prove that
\begin{equation}
\label{eqn:**}
a\,T_{\mathbf i^c\boldsymbol\ell^c}
 =T_{\mathbf i\boldsymbol\ell}.
\end{equation}
Let $\mathbf i, \boldsymbol \ell \subset \{1,\ldots,4n\}$ be increasing $2n$-tuples, and consider the disjoint sets
$$
C=\mathbf i\cap\boldsymbol\ell\comma
A=\mathbf i\cap\boldsymbol\ell^c\comma
B=\mathbf i^c\cap\boldsymbol\ell\comma
D=\mathbf i^c\cap\boldsymbol\ell^c.
$$
Then $s\defeq |C|=|D|$, while $|A|=|B|=2n-s$. If $s=0$, then
$\boldsymbol\ell=\mathbf i^c$ and
$\boldsymbol\ell^c=\mathbf i$, so \eqref{eqn:**} follows immediately from block
symmetry. Suppose therefore that $s\geq1$, set $E=C\cup D$, and, for
each $s$-element subset $P\subset E$, define
$$
x_{\scalebox{0.5}{$P$}}\defeq T_{AP,BP},
$$
where juxtaposition denotes concatenation of ordered tuples; note that $i \in E$ if and only if $i\notin A$ and $i\notin B$. If
$Q\subset E$ has $s-1$ elements, then total trace-freeness gives
\begin{equation}
\label{eqn:T**}
0=(cT)_{AQ,BQ}
  =\sum_{m=1}^{4n}T_{mAQ,mBQ}
  =\sum_{m\in E\setminus Q}x_{\scalebox{0.5}{$Q\cup\{m\}$}}.
\end{equation}
Indeed, the terms with $m\in A$, $m\in B$, or $m\in Q$ vanish by
alternation, while moving $m$ past $A$ and past $B$ produces two equal
signs, which cancel.
For $0\leq r\leq s$, let
$$
X_r\defeq 
\sum_{\substack{P\subset E,\ |P|=s\\|P\cap C|=r}}x_{\scalebox{0.5}{$P$}}.
$$
Summing \eqref{eqn:T**} over all $(s-1)$-element subsets $Q\subset E$ with
$|Q\cap C|=r$ yields
$$
(r+1)X_{r+1}+(s-r)X_r=0\comma 0\leq r\leq s-1.
$$
It follows that $
X_r=(-1)^r\binom{s}{r}X_0$, hence $X_s=(-1)^sX_0$. Since $X_s = x_{\scalebox{0.5}{$C$}}$ and $X_0 = x_{\scalebox{0.5}{$D$}}$ (recall that $C$ and $D$ are disjoint), we thus have that
$
x_{\scalebox{0.5}{$C$}}=(-1)^sx_{\scalebox{0.5}{$D$}}.
$
If $\rho(U,V)$ denotes the sign required to rearrange the
concatenation $UV$ into increasing order, then
$$
x_{\scalebox{0.5}{$C$}}=T_{AC,BC}=\rho(A,C)\rho(B,C)T_{\mathbf i\boldsymbol\ell}\commass
x_{\scalebox{0.5}{$D$}}=T_{AD,BD}=\rho(A,D)\rho(B,D)T_{\mathbf i^c\boldsymbol\ell^c},
$$
where block symmetry was used in the second equality. Next, using that $\vep_{\scalebox{0.5}{$ACBD$}}=\rho(A,C)\rho(B,D)\varepsilon_{\mathbf i\mathbf i^c}$, that $\vep_{\scalebox{0.5}{$BCAD$}}=\rho(B,C)\rho(A,D)\varepsilon_{\boldsymbol \ell\boldsymbol \ell^c}$, and finally that $\vep_{\scalebox{0.5}{$BCAD$}} = (-1)^{s(2n-s)}\vep_{\scalebox{0.5}{$ACBD$}}=(-1)^{s}\vep_{\scalebox{0.5}{$ACBD$}}$, gives
$$
T_{\mathbf i\boldsymbol\ell}=(-1)^s\rho(A,C)\rho(B,C)\rho(A,D)\rho(B,D)T_{\mathbf i^c\boldsymbol\ell^c}
 =
\varepsilon_{\mathbf i\mathbf i^c}
\varepsilon_{\boldsymbol\ell\boldsymbol\ell^c}T_{\mathbf i^c\boldsymbol\ell^c}.
$$
Thus $T_{\mathbf i\boldsymbol\ell}
=aT_{\mathbf i^c\boldsymbol\ell^c}$, proving \eqref{eqn:**}, and consequently
$*\hat T*=\hat T$.
\end{proof}
Cf.~Lemma \ref{lem:hodge-commutation} with \cite[Corollary~4.2]{labbi2005}. When $n=1$, Lemma \ref{lem:hodge-commutation} reduces to the well known fact that a trace-free
algebraic curvature tensor in dimension four commutes with $*$; in particular, that $[\hat{W},*]=0$ in dimension four.

\begin{defn}[Higher self-duality]
\label{def:higher-self-dual}
Let $(M,g)$ be an oriented Riemannian $4n$-manifold.  Denote the $\Lambda^{2n}_{\pm}$-preserving blocks of the Thorpe--Kulkarni operator $\Kht{2n}$ by
$
\Kht{2n}^{\pm}
\defeq\Kht{2n}\big|_{\Lambda^{2n}_{\pm}}.
$
We call $g$ \emph{higher self-dual} if
$\Kht{2n}^{-}=0$, and \emph{higher anti-self-dual} if
$\Kht{2n}^{+}=0$.
\end{defn}
The terminology agrees with the usual four-dimensional convention because
$\Kt{2}=W$.  Reversing the orientation changes $*$ to $-*$ on middle
forms and therefore exchanges the two conditions. The Hodge splitting is useful because the top Pontryagin form is a signed
trace: the star contributes a plus sign on $\Lambda^{2n}_+$ and a minus
sign on $\Lambda^{2n}_-$.  To make this precise, let
$\mathfrak p_n(g)$ denote the Chern--Weil
form representing $p_n(TM)$.  We normalize it so that its integral over a
closed oriented manifold is the Pontryagin number $p_n(M)$.  Recall that
$
d_n=\frac{(2n)!}{2^n n!},
$
the number of unordered pairings of $2n$ objects.

\begin{thm}[Top Pontryagin identity]
\label{thm:pontryagin}
On every oriented Riemannian $4n$-manifold,
\begin{equation}
\label{eq:pontryagin-pointwise}
\mathfrak p_n(g)
=\frac{d_n^2}{(2\pi)^{2n}}
\left(
|\Kht{2n}^{+}|^{2}
-|\Kht{2n}^{-}|^{2}
\right)\dV.
\end{equation}
If $M$ is closed and
$$
\EK^{\pm}(g)
\defeq\int_M|\Kht{2n}^{\pm}|^{2}\dV,
$$
then
\begin{equation}
\label{eq:energy-difference}
\EK^{+}-\EK^{-}
=\gamma_n p_n(M)\comma
\gamma_n=(2\pi)^{2n}
\left(\frac{2^n n!}{(2n)!}\right)^{\!2}\cdot
\end{equation}
Consequently,
\begin{equation}
\label{eq:sharp-bound}
\EK(g)=\gamma_n|p_n(M)|
+2\min\{\EK^{+}(g),
\EK^{-}(g)\}.
\end{equation}
Equality in $\EK\geq\gamma_n|p_n(M)|$ holds precisely when $g$ is higher
self-dual or higher anti-self-dual: the self-dual case is forced when
$p_n(M)>0$, the anti-self-dual case when $p_n(M)<0$, and if $p_n(M)=0$
equality means $\Kt{2n}=0$.
\end{thm}

\begin{proof}
We begin by recalling Thorpe's operator in the Chern--Weil form. Let $\{e_1,\ldots,e_{4n}\}$ be an oriented orthonormal basis and write
$\Omega_i{}^j$ for the curvature two-form matrix.  Chern--Weil theory
defines the total Pontryagin form from the invariant polynomial
$\det(I+\Omega/2\pi)$; see, e.g., \cite{milnorStasheff}.
Because $\Omega$ is skew-symmetric, every
principal determinant of even size is the square of a Pfaffian.  Taking
the degree-$4n$ part of $\det(I+\Omega/2\pi)$ therefore gives the
standard principal-Pfaffian formula
\begin{equation}
\label{eq:pfaffian-sum}
\mathfrak p_n(g)
=\frac{1}{(2\pi)^{2n}}
\sum_{|\mathbf i|=2n}\operatorname{Pf}(\Omega_{\mathbf i})
\wedge\operatorname{Pf}(\Omega_{\mathbf i}).
\end{equation}
Here the sum runs over all increasing $2n$-tuples
$\mathbf i=(i_1<\cdots<i_{2n})$, and $\Omega_{\mathbf i}$ is the
principal $2n$-by-$2n$ submatrix selected by those indices.  By definition,
$$
\operatorname{Pf}(\Omega_{\mathbf i})
=\frac{1}{2^n n!}\sum_{\sigma\in S_{2n}}
\operatorname{sgn}(\sigma)
\Omega_{i_{\sigma(1)}i_{\sigma(2)}}\wedge\cdots\wedge
\Omega_{i_{\sigma(2n-1)}i_{\sigma(2n)}}.
$$
Thus a Pfaffian is a signed sum over the $d_n$ pairings of its indices.
Let $\mathbf j=(j_1<\cdots<j_{2n})$, put
$e^{\mathbf j}=e^{j_1}\wedge\cdots\wedge e^{j_{2n}}$, and use
$\Omega_{ij}(e_k,e_l)=R_{ijkl}$.  Evaluating the displayed wedge product
on $e_{j_1},\ldots,e_{j_{2n}}$ introduces the second alternating sum in
\eqref{eq:Rp-def}.  The factors in the two definitions then give
$$
\operatorname{Pf}(\Omega_{\mathbf i})(e_{\mathbf j})
=d_n\Rt{2n}(e_{\mathbf i},e_{\mathbf j})\comma
\operatorname{Pf}(\Omega_{\mathbf i})
=d_n\sum_{|\mathbf j|=2n}
(\Rt{2n})_{\mathbf i\mathbf j}e^{\mathbf j}.
$$
Consequently,
$$
\sum_{|\mathbf i|=2n}\operatorname{Pf}(\Omega_{\mathbf i})
\wedge\operatorname{Pf}(\Omega_{\mathbf i})
=d_n^2\sum_{|\mathbf i|=2n}
\langle *\Rht{2n}e_{\mathbf i},
\Rht{2n}e_{\mathbf i}\rangle\dV.
$$
The sum on the right is the operator trace, so \eqref{eq:pfaffian-sum}
gives
\begin{equation}
\label{eq:thorpe-p-form}
\mathfrak p_n(g)
=\frac{d_n^2}{(2\pi)^{2n}}
\tr(*\Rht{2n}^{2})\dV.
\end{equation}
This higher-curvature-operator form of $\mathfrak p_n$ is known; see \cite{stehney,bivens}. We now
simplify \eqref{eq:thorpe-p-form} by showing that the metric-trace pieces of $\Rt{2n}$ do not contribute to the
Pontryagin form, so that $\Rht{2n}$ may be replaced by $\Kht{2n}$ (for a more general metric-factor vanishing identity, see \cite[Lemma~3.7]{labbi2}).
Indeed, for any two generalized curvature tensors $A,B\in\mathscr C_{2n}(V)$, put
$$
\mathcal P(A,B)\defeq\tr(*\hat A\hat B).
$$
With the ordered-tuple convention that we used in Lemma
\ref{lem:hodge-commutation}, this is exactly the complete contraction
$$
\mathcal P(A,B)=\frac{1}{((2n)!)^3}
\varepsilon^{\mathbf i\mathbf k}
A_{\mathbf i}{}^{\mathbf j}B_{\mathbf k\mathbf j}.
$$
Now suppose $A=\G H$.  Each term of $\G H$ contains a metric pairing \mbox{between}
one index in the $\mathbf i$-block and one in the $\mathbf j$-block.
For example, the term
$\delta_{i_1}^{j_1}H_{i_2\ldots i_{2n}}{}^{j_2\ldots j_{2n}}$
identifies $j_1$ with $i_1$ in the complete contraction.  The alternating
symbol then alternates
$B_{k_1\ldots k_{2n}i_1j_2\ldots j_{2n}}$ over the $2n+1$ indices
$k_1,\ldots,k_{2n},i_1$.  The generalized algebraic Bianchi identity \eqref{eqn:gen_bian} says
that this alternation is zero.  Every other term in $\G H$ differs only by
a relabeling and a sign, so it vanishes for the same reason.  Thus
$$
\mathcal P(\G H,B)=0
\qquad\text{for every }B\in\mathscr C_{2n}(V).
$$
The same conclusion holds with the two arguments interchanged; indeed,
because $A,B,$ and $*$ are self-adjoint,
$$
\mathcal P(B,A)=\tr(*\hat B\hat A)
=\tr\big((*\hat B\hat A)^t\big)
=\tr(\hat A\hat B*)=\mathcal P(A,B).
$$
Finally, by Lemma \ref{lem:projection}, we have that $\Rt{2n}=\Kt{2n}+\G H$ (recall \eqref{eq:orthogonal-decomp}), so
\begin{equation}
\label{eq:replace-by-K}
\tr(*\Rht{2n}^{2})
=\tr(*\Kht{2n}^{2}).
\end{equation}

With this in hand, we can now use the fact that $[\Kht{2n},*]=0$; i.e., that $\Kht{2n}$ is block diagonal on
$\Lambda^{2n}_{+}
\oplus\Lambda^{2n}_{-}$, hence so is $\Kht{2n}^2$.  Because it is also self-adjoint,
$$
\tr(*\Kht{2n}^{2})=\tr\big((\Kht{2n}^{+})^2\big)
-\tr\big((\Kht{2n}^{-})^2\big)=|\Kht{2n}^{+}|^{2}
-|\Kht{2n}^{-}|^{2}.
$$
Together with \eqref{eq:thorpe-p-form} and \eqref{eq:replace-by-K}, this proves
\eqref{eq:pontryagin-pointwise}.  Integrating and solving for the
difference of the two energies gives \eqref{eq:energy-difference}.  Finally,
for any nonnegative numbers $x,y$, observe that
$$
x+y=|x-y|+2\min\{x,y\}.
$$
Applying this with $x=\EK^{+}$ and
$y=\EK^{-}$, and noting that $\EK=\EK^++\EK^-$, yields
\eqref{eq:sharp-bound} and its equality statement.
\end{proof}

\begin{cor}
\label{cor:selfdual-bachflat}
Every higher self-dual or higher anti-self-dual metric on a closed
$4n$-manifold is an absolute minimizer of $\EK$ and satisfies
$\Bn=0$.
\end{cor}

\begin{proof}
The lower bound in Theorem \ref{thm:pontryagin} depends only on $TM$ and
the orientation, so that $\EK(g)$ is minimized precisely when $\min\{\EK^{+}(g),
\EK^{-}(g)\}=0$.  Any metric for which one of these blocks vanishes attains
it, so its first variation
vanishes in every metric direction.  By \eqref{eq:first-variation-convention} in Theorem \ref{thm:bach-symmetry}, $\Bn=0$.
\end{proof}

For $n=1$, $\Kt{2}=W$, $d_1=1$, and $p_1(M)=3\tau(M)$.  Thus
\eqref{eq:pontryagin-pointwise} becomes the Hirzebruch signature formula
\cite{hirzebruch}:
$$
\tau(M)=\frac{1}{12\pi^2}\int_M
\left(|\hat W^{+}|^{2}
-|\hat W^{-}|^{2}\right)\dV.
$$
For $n>1$, however, the signature is determined by the full Hirzebruch
$L$-polynomial, not by $p_n$ alone (e.g., $L_2=\frac{1}{45}(7p_2-p_1^2)$, etc.).  To our knowledge, \eqref{eq:pontryagin-pointwise} is a new formula specifically for the top Pontryagin form.

\section{Fubini--Study metrics}
\label{sec:fubini-study}

We now show that higher self-duality has nonzero compact examples in every
dimension $4n$.  The Fubini--Study metric is a natural candidate due to its
large isotropy group: at every point its curvature is fixed by $U(2n)$.
We therefore need only determine where $U(2n)$-fixed vectors can occur in
the two Hodge pieces of the Kulkarni module.

\begin{lemma}[Unitary fixed vectors]
\label{lem:unitary}
Let $V=\mathbb{R}^{4n}$ carry a compatible orthogonal complex structure and the corresponding complex orientation. If
$$
T\in\mathcal{K}_{\scalebox{0.6}{$2n$}}\defeq\{K\in\mathcal{C}_{2n}(V):cK=0\}
$$
is invariant under $U(2n)$, then the restriction of $\hat{T}$ to $\Lambda^{2n}_{-}$ vanishes.
\end{lemma}

\begin{proof}
For $n=1$, this is the familiar fact that a $U(2)$-invariant Weyl tensor lies in $\mathcal{W}^{+}$. Thus let $n\geq 2$ and consider $\mathcal{K}_{\scalebox{0.6}{$2n$}}^{\pm}\defeq\{K\in\mathcal{K}_{\scalebox{0.6}{$2n$}}:*\hat{K}=\pm\hat{K}\}$. The complexification of $\mathcal{K}_{\scalebox{0.6}{$2n$}}^{-}$ is the irreducible $SO(4n)$-module with highest weight
$$
4\varpi_{2n-1};
$$
see \cite[Sections 10.2.3--10.2.5 and Appendix~F]{goodmanWallach}. The Cartan--Helgason theorem for the compact symmetric pair
$
SO(4n)/U(2n)
$
states that the spherical highest weights are generated by
$$
\varpi_{2},\varpi_{4},\ldots,\varpi_{2n-2},2\varpi_{2n};
$$
see \cite[Chapter~V, Theorem~4.1]{helgason} and the explicit type-DIII list in \cite[Lemma~11]{matassa}. In particular, $4\varpi_{2n-1}$ is not spherical. Hence
\beqa
\label{eqn:module}
\bigl(\mathcal{K}_{2n}^{-}\bigr)^{U(2n)}=\{0\}.
\eeqa
Write $\hat{T} = \hat{T}^++\hat{T}^-$, where $\hat{T}^{\pm}\defeq \frac{1}{2}(\hat{T}\pm *\hat{T})$, and let $T^{\pm}$ denote the corresponding tensors. Under the $SO(4n)$-equivariant identification \mbox{$K \mapsto \hat{K}$,} the Hodge block projections preserve $\mathcal{K}_{\scalebox{0.6}{$2n$}}$, so $T^{\pm}\in\mathcal{K}_{\scalebox{0.6}{$2n$}}^{\pm}$. Moreover, both $\hat{T}^{\pm}$ are $U(2n)$-invariant, since $*$ commutes with the action of $U(2n)$. Equation \eqref{eqn:module} therefore gives $T^-=0$, and hence $\hat{T}\big|_{\Lambda^{2n}_{-}} = \hat{T}^{-}=0$.
\end{proof}

\begin{thm}[Fubini--Study self-duality]
\label{thm:FS}
For every $n\geq1$, the Fubini--Study metric on $\CP^{2n}$, with its
complex orientation, satisfies
$$
\Kht{2n}^{-}=0\comma
\Kt{2n}\neq0.
$$
It is therefore a nontrivial higher self-dual metric and an absolute
minimizer of $\EK$.  Reversing the orientation makes it higher
anti-self-dual.
\end{thm}

\begin{proof}
At any point of $\CP^{2n}$, the isotropy group of the Fubini--Study metric
is $U(2n)$.  Its curvature tensor is $U(2n)$-invariant, and the operations
used to form $\Rt{2n}$ and then remove its traces commute with the action
of $U(2n)$.  Consequently, $\Kt{2n}$ is $U(2n)$-invariant; by Lemma
\ref{lem:unitary}, 
$\Kht{2n}^{-}=0$. That $\Kt{2n}\neq0$ follows immediately from Theorem \ref{thm:pontryagin} and the fact that $p_n(\CP^{2n})=\binom{2n+1}{n}\neq0$. The minimization statement is
Corollary \ref{cor:selfdual-bachflat}, and orientation reversal exchanges
the two Hodge summands.
\end{proof}

\section{The Einstein condition}
\label{sec:products}
The following compact example, involving round metrics $\mathbb{S}^{2m}(k)$ of curvature $k$, shows that Einstein metrics need not be higher Bach-flat.

\begin{prop}
\label{prop:einstein-counterexample}
The product
$\mathbb{S}^2(3)\times \mathbb{S}^2(3)\times \mathbb{S}^4(1)
$
is Einstein but not higher Bach-flat.
\end{prop}

\begin{proof}
We compute $\Kht4$ first on $\mathbb{S}^2(a)\times \mathbb{S}^2(b)\times \mathbb{S}^4(c)$.
At a point $x$, write
$$
T_xM=V_1\oplus V_2\oplus V_3\comma
(\dim V_1,\dim V_2,\dim V_3)=(2,2,4),
$$
and choose an orthonormal product basis $\{e_1,\dots,e_8\}$.  Let us say that \mbox{$e_{i_1} \wedge \cdots \wedge e_{i_p}$} has
``type $(r,s,u)$" when it uses $r,s,u$ vectors from the three factors;
its multiplicity is $\binom2r\binom2s\binom4u$.  The product isotropy
group preserves $\Rt4$ and $\Kt4$.  Coordinate reflections force $\Rht{4}$ and $\Kht{4}$ to be diagonal, and rotations within the factors
show that each $e_i \wedge e_j \wedge e_k \wedge e_{l}\neq 0$ is an eigenvector whose eigenvalue depends only on $(r,s,u)$.  Then \eqref{eq:R4-explicit} yields
$$
\Rt4(e_i,e_j,e_k,e_l;e_i,e_j,e_k,e_l)
=\frac13(k_{ij}k_{kl}+k_{ik}k_{jl}+k_{il}k_{jk}),
$$
where $k_{ij}$ is the sectional curvature of the two-plane $\ww{e_i}{e_j}$.
Since mixed planes have curvature zero, the only nonzero eigenvalues of
$\Rht4$ are
\beqa
\label{eqn:eig}
\frac{3c^2}{3}\text{ on }(0,0,4)\commas
\frac{ab}{3}\text{ on }(2,2,0)\commas
\frac{bc}{3}\text{ on }(0,2,2)\commas
\frac{ac}{3}\text{ on }(2,0,2).
\eeqa
We now wish to determine the $\Kht4$-eigenvalues $\kappa_{rsu}$ on type
$(r,s,u)$. Consider a three-form $\xi\defeq e_i\wedge e_j \wedge e_k$ of
type $(r,s,u)$, and the contraction $(c\Kt{4})(\xi,\xi) = \sum_{\alpha=1}^8 \Kt{4}(e_\alpha,e_i,e_j,e_k,e_\alpha,e_i,e_j,e_k)$.  Since
$c\Kt{4}=0$, this gives
$$
(2-r)\kappa_{r+1,s,u}
+(2-s)\kappa_{r,s+1,u}
+(4-u)\kappa_{r,s,u+1}=0.
$$
For the eight possible three-form types, these equations are, with each $\xi$-type indicated,
$$
\begin{alignedat}{2}
(0,0,3)\colon\hspace{.05in}2\kappa_{103}+2\kappa_{013}+\kappa_{004}&=0,&\qquad
(0,1,2)\colon\hspace{.05in}2\kappa_{112}+\kappa_{022}+2\kappa_{013}&=0,\\
(0,2,1)\colon2\kappa_{121}+3\kappa_{022}&=0,&
(1,0,2)\colon\kappa_{202}+2\kappa_{112}+2\kappa_{103}&=0,\\
(1,1,1)\colon\kappa_{211}+\kappa_{121}+3\kappa_{112}&=0,&
(1,2,0)\colon\kappa_{220}+4\kappa_{121}&=0,\\
(2,0,1)\colon2\kappa_{211}+3\kappa_{202}&=0,&
(2,1,0)\colon\kappa_{220}+4\kappa_{211}&=0.
\end{alignedat}
$$
All of these can be expressed in terms of the single eigenvalue $\kappa_{022}$:
\beqa
\label{eqn:kappa}
\kappa_{013}=\kappa_{103}=\kappa_{121}=\kappa_{211}=-\frac32\kappa_{022}\comma \left.\begin{array}{lr}
\kappa_{112}=\kappa_{202}=\kappa_{022},\\
\kappa_{004}=\kappa_{220}=6\kappa_{022}.
\end{array}\right.
\eeqa
If $V\defeq V_1\oplus V_2 \oplus V_3$ and $G \defeq O(V_1)\times O(V_2) \times O(V_3)$ is the orthogonal group, then by \eqref{eqn:kappa} the space \mbox{$\{T\in\mathscr C_{4}(V) : cT=0\ ,\ g\cdot T = T~\text{for all $g \in G$}\}$} is at most one-dimensional (note that the eigenvectors are \mbox{$e_i \wedge e_j \wedge e_k \wedge e_{l}\neq 0$,} by product-isotropy invariance).  This space is indeed one-dimensional, as can be seen by defining
\mbox{$\mathcal T \defeq \sum_{|I|=4}\lambda_{\scalebox{0.6}{$\text{type}(I)$}}e_I\otimes e_I$} by $\lambda_{\scalebox{0.6}{$(0,0,4)$}}=1$; by \eqref{eqn:kappa}, 
its three eigenvalue groups are $-1/4,1/6,1$.  Such a $\mathcal{T}$ is $G$-invariant, satisfies the generalized
algebraic Bianchi identity \eqref{eqn:gen_bian}, and
trace-free by the eight relations above.  Using $\binom2r\binom2s\binom4u$, the multiplicities of
the three eigenvalue groups are, respectively, $32,36,2$ (e.g., types $(1,1,2), (2,0,2)$, and $(0,2,2)$ comprise the 1/6-eigenvalue group, and have multiplicities $24,6,6$, for a total of $36$). Hence
$$
|\hat{\mathcal T}|^2
=32\Big(\frac14\Big)^{\!2}+36\Big(\frac16\Big)^{\!2}+2=5.
$$
Since $K_4$ is the orthogonal projection of $R_4$ onto the trace-free
space, we have that $K_4=\kappa_{004}\mathcal T$ and thus, using \eqref{eqn:eig},
$$
\kappa_{004}|\hat{\mathcal T}|^2
=\langle \Rht{4},\hat{\mathcal T}\rangle
=\Big(6\frac{bc}{3}+6\frac{ac}{3}\Big)\frac16+\Big(\frac{3c^2}{3}+\frac{ab}{3}\Big)
=\frac{ab+ac+bc+3c^2}{3}\cdot
$$
Setting $Q\defeq ab+ac+bc+3c^2$, we obtain
$
\kappa_{004}=Q/15.
$
Hence the $\Kht4$-eigenvalues are $-Q/60$ on
$(0,1,3),(1,0,3),(1,2,1),(2,1,1)$, $Q/90$ on
$(0,2,2),(1,1,2),(2,0,2)$, and $Q/15$ on $(0,0,4),(2,2,0)$, with respective multiplicities
$32,36,2$.  In particular,
\begin{equation}
\label{eq:triple-K4-norm}
|\Kht4|^2
=32\Big(\frac{Q}{60}\Big)^{\!2}
+36\Big(\frac{Q}{90}\Big)^{\!2}
+2\Big(\frac{Q}{15}\Big)^{\!2}
=\frac{Q^2}{45}\cdot
\end{equation}

Now take $a=b=3$ and $c=1$.  Since
$\Ric_{\scalebox{0.6}{$\mathbb{S}^m(k)$}}=(m-1)k\,g$, all three factors have Ricci tensor $3g$;
their product is therefore Einstein with $\Ric_{\scalebox{0.6}{$g$}}=3g$.  Consider the family of constant-curvature metrics obtained by rescaling only
the first two-sphere:
$$
g_t\defeq t g_1\oplus g_2\oplus g_3.
$$
The curvatures are $3/t,3,1$, so
$Q(t)=6+12/t$. Since $dV_{\scalebox{0.6}{$g_t$}}=t^{2/2}\dVg=t\dVg$ and $|\Kht4(t)|_{\scalebox{0.6}{$g_t$}}^2= \frac{Q^2(t)}{45}$ by \eqref{eq:triple-K4-norm} (which quantity is constant over $M$),
$$
\EK(g_t)
=\frac{t}{45}\Big(6+\frac{12}{t}\Big)^{\!2}\Vol(M,g)\Rightarrow
\left.\frac{d}{dt}\right|_{t=1}\EK(g_t)
=-\frac{12}{5}\Vol(M,g)\neq0.
$$
Thus $g$ is not a critical point of $\EK$, hence \eqref{eq:first-variation-convention} yields $\mathcal{B}^{\scalebox{0.5}{$(2)$}}\neq0$.
\end{proof}

Kastor proved that a metric is higher Bach-flat when it is
``Einstein of order $n$" \cite[Eq.~(39)]{kastor}.  In our notation, this says that the first contraction of $\Rt{2n}$, viewed as an operator, is scalar: $\widehat{cR}_{\scalebox{0.6}{$2n$}}=\lambda\id_{\Lambda^{2n-1}}$. For $n=1$ this is the ordinary Einstein condition (recall \eqref{eqn:Ricci}), but not for $n > 1$.

\section*{References}
\renewcommand*{\bibfont}{\footnotesize}
\printbibliography[heading=none]
\end{document}